\documentclass[12pt]{amsart}
\usepackage{amsmath, amssymb, amsthm,
amsfonts, amsthm, bm, eucal, graphicx}
\usepackage{xypic}
\xyoption{all}

\theoremstyle{plain}
\newtheorem{thm}{Theorem}[section]
\newtheorem{cor}[thm]{Corollary}
\newtheorem{lem}[thm]{Lemma}
\newtheorem{prop}[thm]{Proposition}

\theoremstyle{definition}
\newtheorem{defi}[thm]{Definition}
\newtheorem{conj}[thm]{Conjecture}
\newtheorem{nota}[thm]{Notation}
\newtheorem{obs}[thm]{Observation}
\newtheorem{obss}[thm]{Observations}
\newtheorem{rem}[thm]{Remark}
\newtheorem{rems}[thm]{Remarks}
\newtheorem{exa}[thm]{Example}
\newtheorem{exas}[thm]{Examples}
\newtheorem{prob}[thm]{Problem}
\newtheorem{sit}[]{}

\newcommand{\brem}{\begin{rem}}
\newcommand{\brems}{\begin{rems}}
\newcommand{\erem}{\end{rem}}
\newcommand{\erems}{\end{rems}}
\newcommand{\bexa}{\begin{exa}}
\newcommand{\bexas}{\begin{exas}}
\newcommand{\eexa}{\end{exa}}
\newcommand{\eexas}{\end{exas}}
\newcommand{\bdefi}{\begin{defi}}
\newcommand{\edefi}{\end{defi}}
\newcommand{\bdefis}{\begin{defis}}
\newcommand{\edefis}{\end{defis}}
\newcommand{\bcor}{\begin{cor}}
\newcommand{\ecor}{\end{cor}}
\newcommand{\blem}{\begin{lem}}
\newcommand{\elem}{\end{lem}}
\newcommand{\bconv}{\begin{conv}}
\newcommand{\econv}{\end{conv}}
\newcommand{\bconj}{\begin{conj}}
\newcommand{\econj}{\end{conj}}
\newcommand{\bprop}{\begin{prop}}
\newcommand{\eprop}{\end{prop}}
\newcommand{\bprob}{\begin{prob}}
\newcommand{\eprob}{\end{prob}}
\newcommand{\bthm}{\begin{thm}}
\newcommand{\ethm}{\end{thm}}
\newcommand{\bnota}{\begin{nota}}
\newcommand{\enota}{\end{nota}}
\newcommand{\bobs}{\begin{obs}}
\newcommand{\eobs}{\end{obs}}
\newcommand{\bobss}{\begin{obss}}
\newcommand{\eobss}{\end{obss}}
\newcommand{\bsit}{\begin{sit}}
\newcommand{\esit}{\end{sit}}
\newcommand{\be}{\begin{eqnarray}}
\newcommand{\ee}{\end{eqnarray}}
\newcommand{\bproof}{\begin{proof}}
\newcommand{\eproof}{\end{proof}}

\def\ba{\begin{array}}
\def\ea{\end{array}}

\newcommand{\F}{{\mathbb F}}

\newcommand{\N}{{\mathbb N}}
\newcommand{\NO}{{\mathbb N}_{\rm odd}}

\newcommand{\PP}{{\mathbb P}}
\newcommand{\T}{{\mathbb T}}

\newcommand{\C}{{\mathbb C}}
\newcommand{\Q}{{\mathbb Q}}
\newcommand{\Z}{{\mathbb Z}}
\newcommand{\cF}{{\mathcal F}}
\renewcommand{\phi}{\varphi}

\newcommand{\G}{{\Gamma}}
\newcommand{\D}{{\Delta}}
\newcommand{\ord}{{\operatorname{ord}}}
\newcommand{\Harm}{{\operatorname{Harm}}}
\newcommand{\sord}{{\operatorname{subord}}}
\newcommand{\ind}{{\operatorname{ind}}}
\newcommand{\Char}{{\operatorname{Char}}}
\newcommand{\mult}{{\operatorname{mult}}}
\newcommand{\Span}{{\operatorname{span}}}
\newcommand{\card}{{\operatorname{card}}}
\newcommand{\red}{{\operatorname{red}}}

\newcommand{\vlin}[1]{\hspace{0.75mm}\unitlength1mm\begin{picture}
(#1,0)
                     \put(0,0){\line(1,0){#1}}
                    \end{picture}\hspace{0.75mm}\rule[-3mm]{0mm}{4mm}}

\def\llin{\vlin{11.5}}

\newcommand{\co}[1]{\unitlength1mm\begin{picture}(0,8)
  \put(0,0){\circle{1.5}}
  \put(0,3){\makebox(0,5)[b]{$#1$}}
                    \end{picture}}
\newcommand{\mybox}{\unitlength1mm\begin{picture}(0,1.5)
  \put(-0.75,-0.75){\line(0,1){1.5}}
  \put(-0.75,-0.75){\line(1,0){1.5}}
  \put(0.75,0.75){\line(0,-1){1.5}}
  \put(0.75,0.75){\line(-1,0){1.5}}
  \end{picture}}
\newcommand{\boxo}[1]{\unitlength1mm\begin{picture}(0,8)
  \put(0,0){\mybox}
  \put(0,3){\makebox(0,5)[b]{$#1$}}
                    \end{picture}}

\title[Periodic harmonic functions on lattices in positive
characteristic]
{Periodic harmonic functions on lattices\\
and points count in positive characteristic}

\author{Mikhail Zaidenberg}
\address{Universit\'e
Grenoble I, Institut Fourier, UMR 5582 CNRS-UJF, BP 74, 38402 St.\
Martin d'H\`eres c\'edex, France} \email{zaidenbe@ujf-grenoble.fr}

\thanks{
\mbox{\hspace{11pt}}{\it 1991 Mathematics Subject
Classification}: 11B39, 11T06, 11T99, 31C05, 37B15, 43A99.\\
\mbox{\hspace{11pt}}{\it Key words}: cellular automaton,
Chebyshev-Dickson polynomial, convolution operator, lattice,
finite field, discrete Fourier transform, discrete harmonic
function, pluri-periodic function.}

\thanks{{\bf Acknowledgements:} This survey is based on 2 author's
preprints at the Max-Planck-Institute of Mathematics (Bonn), and
on a talk on the meeting ``Analysis on Graphs and Fractals", the
Cardiff University, 29 May-2 June 2007 (a satellite meeting of the
programme ``Analysis on Graphs and its Applications" at the Isaac
Newton Institute from 8 January to 29 June 2007). The author
thanks all these institutions for a generous support.}

\begin{document}

\begin{abstract} This survey addresses
pluri-periodic harmonic functions on lattices with values in a
positive characteristic field.  We mention, as a motivation, the
game ``Lights Out" following the work of Sutner \cite{Su},
Goldwasser-Klostermeyer-Ware \cite{GKW}, Barua-Ramakrishnan-Sarkar
\cite{BR, SB}, Hunzikel-Machiavello-Park \cite{HMP} e.a.; see also
\cite{Za1, Za2} for a more detailed account. Our approach explores
harmonic analysis and algebraic geometry over a positive
characteristic field. The Fourier transform allows us to interpret
pluri-periods of harmonic functions on lattices as torsion
multi-orders of points on the corresponding affine algebraic
variety.
\end{abstract}

\maketitle
\date{}

\section{Introduction}
We consider the Caley graph $\Gamma$ of a free abelian group
(i.e., a lattice) and harmonic functions on $\Gamma$ with values
in a field $K$ of positive characteristic. We are interested in
determining all pluri-periods of such functions. In the
characteristic 2 case, this question arises naturally in relation
with the game ``Lights Out" on a rectangular or a toric board, or
otherwise in studies on the dynamics of linear cellular automata
on a lattice $\Lambda$. We present two possible reductions of this
problem. The first one, developed by Sutner \cite{Su},
Goldwasser-Klostermeyer-Ware \cite{GKW}, Barua-Ramakrishnan-Sarkar
\cite{BR, SB}, Hunzikel-Machiavello-Park \cite{HMP} e.a. deals
with the Chebyshev-Dickson polynomials and their generalizations.
The second one leads to points count on a certain affine algebraic
variety $\Sigma$ over the algebraic closure of $K$. The points on
$\Sigma$ correspond to harmonic characters on $\Lambda$. We
express the pluri-periods of harmonic functions on $\Lambda$, or,
which is the same, the sizes of the toric boards obstructed for
the ``Lights Out" game, as torsion multi-orders of the
corresponding points on $\Sigma$.

It is our pleasure to thank Don Zagier for useful suggestions.

\subsection{The game ``Lights Out"}

\smallskip

\noindent The ``Lights Out" is a solitary game on a rectangular
$m\times n$ board. Initially the board is filled in with $0$'s
(``black") and $1$'s (``white"). The rule of the game consists in
the following:
 {\it A click in a cell changes the state to the opposite
in this cell and in all its horizontal and vertical neighbors.}
The goal of the game is: {\it to reach finally the ``all white"
pattern.}

\smallskip

 As an example, let us consider a $3\times 3$ board.
Starting with an initial pattern as shown below and performing a
sequence of clicks in the cells indicated over arrows, we obtain:

$$\begin{pmatrix}
  0 & 0 & 0 \\
  1 & 1 & 0 \\
  0 & 0 & 1
\end{pmatrix}\,\quad \stackrel{(1,3)}{\rightsquigarrow\rightsquigarrow}\quad \begin{pmatrix}
0 & 1 & 1 \\
1 & 1 & 1 \\
0 & 0 & 1
\end{pmatrix}\quad \stackrel{(1,2)}{\rightsquigarrow\rightsquigarrow}\quad
\begin{pmatrix}
1 & 0 & 0 \\
1 & 0 & 1 \\
0 & 0 & 1
\end{pmatrix}\,$$

$$\quad \stackrel{(2,1), (2,3)}{\rightsquigarrow\rightsquigarrow\rightsquigarrow\rightsquigarrow}\quad \begin{pmatrix}
  0 & 0 & 1 \\
  0 & 0 & 0 \\
  1 & 0 & 0
\end{pmatrix}\,\quad \stackrel{(1,3), (3,1)}{\rightsquigarrow\rightsquigarrow\rightsquigarrow\rightsquigarrow}\quad \begin{pmatrix}
0 & 1 & 0 \\
1 & 0 & 1 \\
0 & 1 & 0
\end{pmatrix}\quad \stackrel{(2,2)}{\rightsquigarrow\rightsquigarrow}\quad
\begin{pmatrix}
0 & 0 & 0 \\
0 & 1 & 0 \\
0 & 0 & 0
\end{pmatrix}\,$$
$$\quad \stackrel{(1,2),(2,1),(2,3),(3,2)}{\rightsquigarrow\rightsquigarrow\rightsquigarrow\rightsquigarrow\rightsquigarrow\rightsquigarrow\rightsquigarrow}\quad
\begin{pmatrix}
  0 & 1 & 0 \\
  1 & 1 & 1 \\
  0 & 1 & 0
\end{pmatrix}\,\quad \stackrel{(2,2)}{\rightsquigarrow\rightsquigarrow}\quad \begin{pmatrix}
0 & 0 & 0 \\
0 & 0 & 0 \\
0 & 0 & 0
\end{pmatrix}\quad. $$

\noindent  To find a shorter way it is enough just to cancel the
clicks that are done twice:

$$\begin{pmatrix}
  0 & 0 & 0 \\
  1 & 1 & 0 \\
  0 & 0 & 1
\end{pmatrix}\,\quad \stackrel{(3,2)}{\rightsquigarrow\rightsquigarrow}\quad \begin{pmatrix}
0 & 0 & 0 \\
1 & 0 & 0 \\
1 & 1 & 0
\end{pmatrix}\quad \stackrel{(3,1)}{\rightsquigarrow\rightsquigarrow}\quad
\begin{pmatrix}
0 & 0 & 0 \\
0 & 0 & 0 \\
0 & 0 & 0
\end{pmatrix}\quad. $$
Indeed, the clicks represent commuting involutions.

\bobs\label{o2} More generally, one can play the ``Lights Out"
game on any finite graph $\Gamma$. Sutner's Garden-of-Eden Theorem
\cite{Su} says that {\it starting with the ``all black" pattern on
$\G$ one can always reach the ``all white" pattern.} That is, the
``all black" pattern is winning for any graph $\G$. \eobs

 A general question is: {\it For which graphs $\G$ one
can win the ``Lights Out" game on $\G$ starting with an arbitrary
initial pattern?} In the latter case we say that $\G$ is winning.
This turns out to be equivalent to a spectral problem for the
corresponding Laplacian $\Delta_\G$ on $\G$. Indeed, the nonzero
harmonic functions on
 $\Gamma$
provide obstructions for the ``Lights Out" game on $\Gamma$
to  always win.

\medskip

\bdefi\label{d1} Let $\Gamma$ be a graph and $K$ be an abelian
group. A function  $h$ on the set of vertices of $\G$ with values
in $K$ is called harmonic\footnote{Alternatively, one can define
harmonic functions by the identity
$$h(v)=\sum_{[v,v']\in\Gamma} h(v')\qquad\forall v\in\Gamma\,.$$
The class of harmonic functions remains the same if $K$ is a field
of characteristic $2$ and changes in case $\Char (K)=p>2$.
However, similar results hold after this replacement. We give
below a general approach covering the both cases.} if
$$(\D_\G h)(v):=h(v)+\sum_{[v,v']\in\Gamma} h(v')=0
\qquad\forall v\in\Gamma\,.$$ Or, in other words, if
\be\label{starf} <h,a_v>=0\qquad\forall v\in\Gamma,\qquad\mbox{
where}\qquad a_v=\delta_v +\sum_{[v,v']\in\Gamma}
\delta_{v'}\,.\ee We call $a_v$ the star-function centered at
$v$.\edefi

\bobs\label{o1} Given a finite field $K$ and a simple finite graph
$\Gamma$, one can equally play the game ``Lights Out" on $\Gamma$
with patterns taking values in $K$; see e.g. \cite{GMT}. The click
in a vertex $v\in\Gamma$ corresponds to the translation
$f\longmapsto f+a_v\,$ in the vector space $\cF(\G,K)$ of all
$K$-valued functions on $\Gamma$. Thus $f$ is
winning\footnote{That is, starting with $f$ one can reach the
``all white" pattern.} if and only if $f\in{\rm
span}\,(a_v\,|\,v\in\Gamma)\,.$

For any harmonic function $h$ on $\Gamma$ and for any $f\in
\cF(\G,K)$, by virtue of (\ref{starf})
$$<h,f>\,=\,<h,f+a_v>\qquad\forall v\in\Gamma.$$
Hence $h$ provides a linear invariant of the game ``Lights Out" on
$\Gamma$, and any such invariant appears in this way. Therefore a
pattern $f: \Gamma\to K$ is winning  if and only if $f\,\bot\,{\rm
Harm}(\Gamma, K)\,,$ where ${\rm Harm}(\Gamma, K)=\ker (\D_\G)$
stands for the space of all $K$-valued harmonic functions on
$\Gamma$. So $\Gamma$ is winning if and only if ${\rm
Harm}(\Gamma,K)=(0)\,.$ \eobs

\bprob\label{pr1} Given a finite field $K$, determine all winning
$m\times n$ boards, or, alternatively, all those which possess a
nonzero $K$-valued harmonic function.\eprob

\bexa\label{e3} For $K=\F_2$, the square board $3\times 3$ is
winning, whereas the boards $4\times 4$ and $5\times 5$ are not as
both of them possess nonzero binary harmonic functions, for
instance
$$\qquad
\begin{pmatrix}
1 & 1 & 1 & 0\\
0 & 1 & 0 & 1\\
0 & 0 & 1 & 1\\
0 & 0 & 0 & 1
\end{pmatrix}\qquad\mbox{resp.,}\qquad
\begin{pmatrix}
1 & 1 & 0 & 1 & 1\\
0 & 0 & 0 & 0 & 0\\
1 & 1 & 0 & 1 & 1\\
0 & 0 & 0 & 0 & 0\\
1 & 1 & 0 & 1 & 1
\end{pmatrix}\,\quad.$$\eexa

\bobs\label{o1bis} Problem \ref{pr1} for a rectangular $m\times n$
board is closely related to a similar question for the toric
$(m+1)\times (n+1)$ board; see \cite{Za1}. Considering the game
``Lights Out" on toric boards rather that on rectangular ones
provides certain advantages. Indeed, the toric $m\times n$ board
$\T_{m,n}$ represents the Caley graph of the abelian group $\Z/m\Z
\times \Z/n\Z\,$ for the standard choice of generators. Its
maximal abelian cover is the Caley graph of the free abelian group
$\Lambda=\Z^2$. Every harmonic function $h\in\Harm(\T_{m,n}, K)$
can be lifted to a bi-periodic harmonic function $\tilde
h\in\Harm(\Lambda, K)$ with periods $me_1$ and $ne_2$.  Thus
Problem \ref{pr1} for toric boards is equivalent to the following
one. \eobs

\bprob\label{pr2} Given a lattice $\Lambda$ and a field $K$,
determine  the
pluri-periods of all nonzero pluri-periodic harmonic functions
$h:\Lambda\to K$.
\eprob

\bexa\label{e10} The game ``Lights Out" played over the binary
field $\F_2$ on the toric board $\T_{10,10}$ does not always win.
Indeed, $\T_{10,10}$ possesses nonzero binary harmonic functions,
for instance, the following one obtained via the doubling of
periods trick \cite{Za1}:
$$h=\begin{pmatrix}
  0 & {\bf 1} & {\bf 1} & 0 & {\bf 1} & 0 & {\bf 1} & 0 & {\bf 1}
  & {\bf 1} \\
  {\bf 1} & 0 & 0 & 0 & 0 & 0 & {\bf 1} & 0 & 0 & 0 \\
  {\bf 1} & 0 & {\bf 1} & 0 & {\bf 1} & {\bf 1} & 0 & {\bf 1} & {\bf 1} & 0  \\
  0 & 0 & {\bf 1} & 0 & 0 & 0 & {\bf 1} & 0 & 0 & 0  \\
  {\bf 1} & {\bf 1} & 0 & {\bf 1} & {\bf 1} & 0 & {\bf 1} & 0 & {\bf 1} & 0  \\
  0 & 0 & {\bf 1} & 0 & 0 & 0 & 0 & 0 & {\bf 1} & 0 \\
  {\bf 1} & 0 & {\bf 1} & 0 & {\bf 1} & 0 & {\bf 1} & {\bf 1} & 0 & {\bf 1}  \\
  0 & 0 & 0 & 0 & {\bf 1} & 0 & 0 & 0 & {\bf 1} & 0 \\
  {\bf 1} & 0 & {\bf 1} & {\bf 1} & 0 & {\bf 1} & {\bf 1} & 0 & {\bf 1} & 0  \\
  {\bf 1} & 0 & 0 & 0 & {\bf 1} & 0 & 0 & 0 & 0 & 0
\end{pmatrix}$$
This pattern $h$ composed of five crosses
$$\begin{pmatrix}
  &  & {\bf 1} & &   \\
  &  & {\bf 1} & & \\
  {\bf 1} & {\bf 1} & 0 & {\bf 1} & {\bf 1}   \\
  &  & {\bf 1} & &   \\
  &  & {\bf 1} & &
\end{pmatrix}\,$$
lifts to a bi-periodic binary harmonic function $\tilde h$ on the
lattice $\Lambda=\Z^2$ with periods $10e_1,\,10e_2$.
\eexa

\bobs\label{o4} The Laplacian $\Delta_\Lambda$ acting on the space
$\cF(\Lambda,K)$ of all functions $\Lambda\to K$ provides a linear
cellular automaton on $\Lambda$ \cite{MOW}. Actually any
homogeneous linear cellular automaton on $\Lambda$ appears in this
way. Furthermore, $\Delta_\Lambda$ can be expressed as the
convolution operator
$$\Delta_\Lambda: f\longmapsto f*a_0\,$$ with kernel
the star-function  $a_0$ on $\Lambda$ centered at the origin. Thus
$$\tilde h\in {\rm Harm} (\Lambda,K)
\quad\Longleftrightarrow\quad {\tilde h}*a_0=0\,.$$ The period
vectors of a pluri-periodic harmonic function $\tilde h$ on
$\Lambda$ form a finite index sublattice
$$\Lambda'=\Lambda'(\tilde h)\subseteq \Lambda\,.$$
The quotient $\T=\Lambda/\Lambda'\,$ is a finite abelian group,
and $\tilde h$ is the pull-back of a function $h:\T\to K\,$
harmonic with respect to an appropriate Laplacian $\D_\T$ on $\T$.
\eobs

\bexa\label{e4} For a circular graph $\T_n={\rm Caley}(\Z/n\Z)\,$
and for any field $K$ of characteristic $p>0$, one has
$${\rm Harm}(\T_n,K)\neq (0)\quad\Longleftrightarrow\quad n\equiv 0 \mod
3\,.$$ Indeed if $n\equiv 0 \mod 3$ then
$$h(k):= k (\!\!\!\!\!\!\mod 3)\in K,\qquad k\in\Z/n\Z\,,$$ is a nonzero
$K$-valued harmonic function on $\T_n$. Conversely, we can write
$$\Delta_{\T_n}=1+\tau+\tau^{-1}\,,$$ where
$$\tau: f(x)\longmapsto f(x+1\mod n)\,$$ is the right shift acting on
$\cF(\T_n,K)$. Hence by virtue of the Spectral Mapping Theorem,
$$0\in {\rm spec} (\Delta_{\T_n})\quad\Longleftrightarrow\quad $$
$$\exists \zeta\in \bar K \,:\, \zeta^n=1,\quad
1+\zeta+\zeta^{-1}=0\,$$
 $$\quad\Longleftrightarrow\quad \exists
\zeta\in \bar K \,:\, \zeta\neq 1,\quad\zeta^n=1,\quad
\zeta^{3}=1\,$$
$$\quad\Longleftrightarrow\quad n\equiv 0 \mod 3\,,$$ where
$\bar K$ stands for the algebraic closure of $K$. \eexa

\subsection{``Lights Out" played on circular graphs} For any finite simple
graph $\G$, the matrix of the ``Markov operator" $\Delta_{\G}-1$
in the canonical base of $\delta$-functions $(\delta_v : v\in \G)$
in the lattice $\cF(\G,\Z)$ is just the adjacency matrix of $\G$.
For a circular graph $\T_n$ with $n\ge 3$ vertices and for a
linear path $\PP_n$ with $n\ge 1$ vertices we have, respectively,
$${\rm adj}\,(\T_n)=\left(\begin{array}{ccccccc}
0 & 1 & 0 & \ldots & 0 & 0 & 1\\
1 & 0 & 1 & \ldots & 0 & 0 & 0\\
0 & 1 & 0 & \ddots & 0 & 0 & 0\\
\vdots & \vdots & \ddots  & \ddots & \ddots  & \vdots & \vdots \\
0 & 0 & 0 & \ddots & 0 & 1 & 0\\
0 & 0 & 0 & \ldots & 1 & 0 & 1\\
1 & 0 & 0 & \ldots & 0 & 1 & 0
 \end{array}\right)\,,\quad
 {\rm adj}\,(\PP_n)=\left(\begin{array}{ccccccc}
0 & 1 & 0 & \ldots & 0 & 0 & 0\\
1 & 0 & 1 & \ldots & 0 & 0 & 0\\
0 & 1 & 0 & \ddots & 0 & 0 & 0\\
\vdots & \vdots & \ddots  & \ddots & \ddots  & \vdots & \vdots \\
0 & 0 & 0 & \ddots & 0 & 1 & 0\\
0 & 0 & 0 & \ldots & 1 & 0 & 1\\
0 & 0 & 0 & \ldots & 0 & 1 & 0
 \end{array}\right)\,.$$

 \smallskip

\noindent We let
$$C_n(x):=(-1)^n\det \left({\rm adj}\,(\T_n)-xI_n\right)$$
denote the characteristic
polynomial of ${\rm adj}\,(\T_n)$. Then $C_n=(x-2)D_n \in \Z[x]$,
where
$$D_3=(x+1)^2,\,\,D_4=x^2(x+2)\qquad\mbox{
and}\qquad D_n=xD_{n-1}-D_{n-2}+2\quad\forall n\ge 5\,.$$ In
particular $3\in {\rm spec} \,(\Delta_{\T_n})$ $\forall n\ge 3$,
and the corresponding eigenfunctions are constant functions on
$\T_n$. Furthermore,
$$C_3(-1)=0,\,\,C_4(-1)=C_5(-1)=-3\qquad\mbox{
and}\qquad C_{n+3}(-1)=C_n(-1)\quad\forall n\ge 3\,.$$ Hence over
$\C$, for $n\ge 3$ one has:
$$0\in {\rm spec} (\Delta_{\T_n})\quad\Longleftrightarrow\quad C_n(-1)=0
\quad\Longleftrightarrow\quad n\equiv 0 \mod 3\,.$$ Consequently,
the ``Lights Out" game on the circular graph $\T_{n}$ ($n\ge 3$)
is winning  over a field $K$ of characteristic $p\neq 3$ if and
only if $n\not\equiv 0\mod 3$. While for $p=3$ none of the graphs
$\T_{n}$ ($n\ge 3$) is winning. Every non-winning graph $\T_{n}$
carries a nonzero $K$-valued harmonic function.  For $p=3$ these
are constant functions. For $p\neq 3$ and $n=3k$ the space $\Harm
(\T_{n},K)$ consists of $3$-periodic functions and is spanned by
the function $h$ from Example \ref{e4} above and its shifts.

\subsection{How can one recognize winning boards?}
In Theorem \ref{1} below we mention two different approaches to
Problem \ref{pr1} for the ``Lights Out" game on toric boards. None
of them is explicit. The first one applies over the Galois field
$\F_p$, while the second one deals with its algebraic closure
$\bar \F_p$. See e.g., \cite{Su,GKW,HMP} for the proof of (a) and
\cite{Za1} for (b).

\bthm\label{1} \begin{enumerate}\item[(a)] For a toric graph
$\T_{m,n}={\rm Caley}\,(\Z/m\Z \times \Z/n\Z)\,$ one has
$${\rm Harm}\,(\T_{m,n},\,\F_p)=(0)
\quad\Longleftrightarrow\quad
\gcd\,(C^{(p)}_m(x),C^{(p)}_n(1-x))=1\,.$$  \item[(b)] For a toric
graph $\T_{\bar n}={\rm Caley}\,(\Z/n_1\Z \times\ldots\times
\Z/n_s\Z)\,$, where $\bar n=(n_1,\ldots,n_s)$, one has
$${\rm Harm}\,(\T_{\bar n},\,\F_p)\neq (0)
\quad\Longleftrightarrow\quad \exists\,
(\zeta_1,\ldots,\zeta_s)\in (\bar \F_p^\times)^s\,:\,$$
$$(*)\qquad
1+\sum_{i=1}^s (\zeta_i+\zeta_i^{-1})=0,
\quad\zeta_i^{n_i}=1,\quad
i=1,\ldots,s \,.$$
\end{enumerate}\ethm

\subsection{Generalized Chebyshev-Dickson polynomials}
Consider further a lattice $\Lambda$, a field $K$
 of characteristic
$p>0$, an arbitrary function $a:\Lambda\to K$  with finite
support, and
 the corresponding Laplacian $$\Delta_a : f \to f*a\,.$$
Let $f:\Lambda\to K$ be a pluri-periodic function with the lattice
of periods $\Lambda'\subseteq \Lambda$. Then clearly the period
lattice of the function $\Delta_a (f)$ contains $\Lambda'$. So the
subspace $\mathcal F (\Lambda, K)^{\Lambda'}$ of all
$\Lambda'$-periodic functions on $\Lambda$ is
$\Delta_a$-invariant, of dimension
$$\dim \mathcal F (\Lambda, K)^{\Lambda'}=\ind
(\Lambda',\Lambda)\,.$$ \bdefi\label{d2} We call a {\it
generalized Chebyshev-Dickson polynomial $T_{a,\Lambda'}$} the
characteristic polynomial of the restriction $\Delta_a|\mathcal F
(\Lambda,K)^{\Lambda'}$. It has degree
$$\deg(T_{a,\Lambda'})=\ind (\Lambda',\Lambda)\,.$$
The classical Chebyshev-Dickson polynomials $T_n$ \footnote{See
Appendix below.} correspond to
$$p=2,\quad\Lambda=\Z,\quad \Lambda'=n\Z,\quad\mbox{and}\quad
a=a_0=\delta_0+\delta_1+\delta_{-1}\,.$$
Given a  base $\mathcal V=(v_1,\ldots,v_s)$ of the lattice $\Lambda$
and a product sublattice $\Lambda'\subseteq
\Lambda$, where
$$\Lambda'=\sum_{i=1}^s
n_i\Z v_i\,,$$ the Chebyshev-Dickson polynomial $T_{a,\Lambda'}$
can be expressed via iterated resultants \cite{Za2}.\edefi

\noindent Like in the classical case, the system of generalized
Chebyshev-Dickson polynomials possesses the following divisibility
properties \cite{Za2}.

\bthm\label{3} \begin{enumerate} \item[(a)]
$\quad\Lambda'\subseteq \Lambda''\quad\Longrightarrow\quad
T_{a,\Lambda''}\,|\,T_{a,\Lambda'}\,.$
\item[(b)] $\ind
(\Lambda',\Lambda'')=p^\alpha\quad\Longrightarrow\quad
T_{a,\Lambda'}=(T_{a,\Lambda''})^{p^\alpha}\,.$\end{enumerate}
\ethm

\subsection{The partnership graph}
In this subsection we return to the special case related to the
game ``Lights Out", where $K=\F_2$, $\Lambda=\Z^2$, and $a=a_0$ is
the star function on $\Lambda$. The covering $\T_{km,ln}\to
\T_{m,n}$ yields an inclusion
$${\rm Harm}\,(\T_{m,n})\hookrightarrow {\rm Harm}\,(\T_{km,ln})\,.$$
Thus one can stick in Problem \ref{pr2} to ``primitive"
non-winning $m\times n$ toric boards.

\bdefi\label{parte} A pair $(m,n)\in\N^2$ is called a pair of
partners if there exists a solution $(\zeta_1,\zeta_2)$ of $(*)$
with exact torsion orders
$$m=\ord (\zeta_1)\quad\mbox{and}\quad n=\ord (\zeta_2)\,.$$
Since $\Char (K)=2$, $m$ and $n$ are odd integers. \edefi

Following a suggestion by Don Zagier, we can represent the above
partnership relation on a ``partnership graph". This graph
$\mathcal P$ has the set $\NO$ of all positive odd integers as the
set of vertices and the pairs of partners for the edges. We label
$\mathcal P$ by attributing to an edge $[m,n]$ the number of
solutions of $(*)$ divided by 2. Given a vertex $n\in\NO$, the sum
of labels over all its incident edges\footnote{A loop at a vertex
is count as a single incident edge.} equals $\varphi(n)$, where
$\varphi$ stands for the Euler totient function. Indeed, given a
primitive $m$th root of unity $\zeta_1$, the equation (*) with
$s=2$ admits exactly 2 solutions of the form $(\zeta_1,\zeta_2)$
and $(\zeta_1,\zeta_2^{-1})$, which yields the claim. In
particular, $\mathcal P$ does not possess isolated vertices. The
following simple observation is also due to Don Zagier.

\bprop\label{2}
All connected components of the partnership graph $\mathcal P$ are
finite.
\eprop

\bproof Given $n\in\NO$, the order and the suborder  of 2 modulo
$n$ are, respectively,
$$f(n)=\ord_n\, 2=\min\{j \,:\, 2^j\equiv 1 \!\!\mod  n\}$$ and
$$ f_0(n)=\sord_n\, 2=\min\{j\,:\, 2^j\equiv
\pm 1 \!\!\mod  n\}\,.$$ Thus $f(n)/f_0 (n)\in\{1,2\}$. Furthermore,
$$f(n)= 2 f_0(n)\quad\mbox{is even}\,\,\,\iff \,\,\,
\exists j\in\N\,:\, 2^j\equiv
-1\!\! \mod  n \,.$$ Letting $q=2^{f_0 (n)}$, $n$ divides exactly
one of $q-1$ and $q+1$. Namely $n\mid (q-1)$ if $f_0 (n)=f(n)$ and
$n\mid (q+1)$ otherwise.

According to (*) a pair $(m,n)$ of odd naturals is a pair of
partners (that is, $[m,n]$ is an edge of ${\mathcal P}$) if and
only if $\xi+\xi^{-1}=1+\eta+\eta^{-1}$ for some primitive roots
$\xi\in\mu_m$ and $\eta\in\mu_n$ (cf. also Example \ref{e1}
below). Thus for a pair of partners $(m,n)$,
$$f_0(m)=\deg (\xi+\xi^{-1})=\deg (\eta+\eta^{-1})=f_0 (n)\,.$$
Hence the suborder function $f_0$ is constant on each connected
component of ${\mathcal P}$.

We let ${\mathcal V}_r=f_0^{-1}(r)$, $r=1,2,\ldots$, denote the
level sets of  $f_0$. By definition of $f_0$, a level set
${\mathcal V}_r$ is contained in the set of all divisors of
$2^{2r}-1$. Therefore it is finite.

Given $n\in\NO$, we let ${\mathcal P}(n)$ denote the
connected component of ${\mathcal P}$ which contains the
vertex $n$. Since the set of vertices of ${\mathcal P}(n)$ is
contained in the finite set ${\mathcal V}_r$, where $r=f_0(n)$,
the former set is finite as well, as
stated.
\eproof

The first 12 level sets $V_r=f_0^{-1}(r)$, $r=1,\ldots,12$, and
the corresponding subgraphs of the labelled partnership graph
${\mathcal P}$ are shown on Figures 1-3 below; they were computed
by Don Zagier with PARI.  A vertex $n$ on these figures is
underlined iff $f_0(n)\neq f(n)$. These computations suggest that
among the $V_r$'s, only $V_5$ is disconnected.

\subsection{Symbolic variety} From now on we let
$K$ be an algebraically closed field  of characteristic $p>0$.
Given a base $\mathcal V=(v_1,\ldots,v_s)$ of a lattice $\Lambda$,
one can identify $\Lambda$ with $\Z^s$, where $s= {\rm rk}\,
(\Lambda)\,$. For a function  $a:\Lambda\to K$ with finite
support, the symbol of the  corresponding Laplacian $\Delta_a$ is
the Laurent polynomial $$\sigma_{a}=
\sum_{u=(u_1,\ldots,u_s)\in\Z^s} a(u)x^{-u}\in
K[x_1,x_1^{-1},\ldots, x_s,x_s^{-1}]$$ with the coefficient
function $a$. The symbolic variety associated with $\Delta_a$ is
$$\Sigma_a=\sigma_a^{-1}(0)\,.$$
More generally, to a sequence $\bar a=(a_1,\ldots,a_t)$
\footnote{In other words, to the system of corresponding
Laplacians $\Delta_{a_1},\ldots,\Delta_{a_t}$.} we associate its
symbolic variety
$$\Sigma_{\bar a}=\{\sigma_{a_j}=0\,:\,j=1,\ldots,t\}
\,,$$ which is a closed subvariety of the affine algebraic torus
$(K^\times)^s$.

\newpage

\begin{center}
\includegraphics{hawo1eps.epsi}
\end{center}
\begin{figure}[h]
\centerline{\hbox{}} \caption[]{} \label{Diagram 1}
\end{figure}

$\,$ \vskip 1in

\begin{center}
\includegraphics{hawo2eps.epsi}
\end{center}
\begin{figure}[h]
\centerline{\hbox{}} \caption[]{} \label{Diagram 2}
\end{figure}

$\,$ \vskip 1in

\begin{center}
\includegraphics{hawo3eps.epsi}
\end{center}
\begin{figure}[h]
\centerline{\hbox{}} \caption[]{} \label{Diagram 3}
\end{figure}

\bexa\label{e1} Consider again $K={\bar \F}_2$ and the Laplacian
$\Delta_{a_0}$ on the plane lattice $\Lambda=\Z^2$ with kernel the
star function $a_0$. The corresponding symbolic variety is the
elliptic cubic curve
$$\Sigma_{a_0}=\{x+1/x+y+1/y=1\}\subseteq (K^\times)^2\,.$$
The logarithm of the Hasse-Weil zeta-function counts points on
$\Sigma_{a_0}$ according to the filtration $\bar
F_2=\bigcup_{n\in\N} \F_{2^n}$. This formula suggests that the
number of toric $m\times n$ boards which admit a nonzero binary
harmonic function is infinite. Moreover, the number of primitive
boards (i.e., those which are not produced using smaller ones) is
also infinite. Indeed, the number of edges of the partnership
graph $\mathcal P$ is infinite, because the number of vertices is
and $\mathcal P$ has no isolated vertex. \eexa

We consider the algebraic closure $K={\bar\F}_p$ of a Galois field
$\F_p$. For $n\in\N$ coprime with $p=\Char (K)$ we let
$\mu_n\subseteq K^\times$ denote the subgroup of $n$th roots of
unity. For a multi-index $\bar n\in\N^s$, where $n_i\not\equiv
0\mod p\,\,\forall i$, we consider the finite $s$-torus
$$\mu_{\bar n}=\mu_{n_1}\times \ldots\times \mu_{n_s}\subseteq
(K^\times)^s\,.$$ The multiplicative group $K^\times$ being a
torsion group, the torus $(K^\times)^s$ is filtered by its finite
subgroups:
$$(K^\times)^s=\bigcup_{\bar n} \mu_{\bar n}\,.$$
Furthermore,
$$(K^\times)^s=\coprod_{\bar n}\nu_{\bar n}\,,$$
where $\nu_{\bar n}\subseteq \mu_{\bar n}$ denotes the set of all
elements of $\mu_{\bar n}$ whose $i$th coordinates are primitive
$n_i$th roots of unity, $i=1,\ldots,s$. Given an algebraic
subvariety $\Sigma\subseteq (K^\times)^s$ we wonder whether the
multi-sequence $\card (\Sigma\cap \nu_{\bar n})$ admits a
recursive generating function.

\subsection{Harmonic characters}
We let $\Char (\Lambda,K^\times)$ denote the set of all characters
$\chi:\Lambda\to K^\times$. Given a base $\mathcal
V=(v_1,\ldots,v_s)$ of $\Lambda$ we consider the associated
isomorphism \be\label{j} j:\Char
(\Lambda,K^\times)\stackrel{\cong}{\longrightarrow}
(K^\times)^s,\qquad \chi\longmapsto
(\chi(v_1),\ldots,\chi(v_s))\,.\ee For $K= {\bar F}_p$ every
$K^\times$-valued character of $\Lambda$ is pluri-periodic. Given
a sublattice $\Lambda'\subseteq \Lambda$ of finite index, all
$\Lambda'$-periodic $K^\times$-valued characters can be produced
by pulling back the $K^\times$-valued characters of the quotient
group $\T=\Lambda/\Lambda'$.

A character $\chi$ is called $a$-harmonic if $\Delta_a (\chi)=0$.
The set of all $a$-harmonic characters of $\Lambda$ is denoted by
$\Char_{a-{\rm harm}} (\Lambda,K^\times)$. The next proposition
follows immediately by using the Fourier transform on a finite
abelian group; see \cite{Za2}.

\bprop\label{p1}
For any product sublattice $\Lambda'=\sum_{i=1}^s
n_i\Z v_i\subseteq \Lambda$  of index
$$\ind (\Lambda',\Lambda)\not\equiv 0 \mod p\,,$$ the space ${\rm
Harm}_a\,(\T_{\bar n},K)$ of all $a$-harmonic
functions on the quotient group
$\T_{\bar n}=\Lambda/\Lambda'$ possesses an
orthonormal basis of $a$-harmonic characters. In particular
$${\Harm}_a\,( \T_{\bar n},K)=\Span\,
\left(\Char_{a-{\rm harm}}\,( \T_{\bar n},K^\times)\right)$$ and
so
$${\Harm}_a\,( \T_{\bar n},K)\neq (0)\quad\Longleftrightarrow\quad
\Char_{a-{\rm harm}}\,( \T_{\bar n},K^\times)\neq\emptyset\,.$$
\eprop

\brem\label{r2}  The latter conclusion remains valid for any (not
necessarily algebraically closed) field $K$ of positive
characteristic. Indeed, the space ${\rm Harm}_a\,( \T_{\bar n},
K)$ is spanned by the traces of harmonic characters; see
\cite{Za2}. \erem

There is a natural bijection between the $a$-harmonic characters
$\Char_{a-{\rm harm}} (\T_{\bar n},K^\times)$ and the points on
the corresponding symbolic variety $\Sigma_a$ with torsion
multi-order dividing $\bar n$. More precisely, the following hold
\cite{Za2}.

\bthm\label{4} Consider a product sublattice
$$\Lambda'=\sum_{i=1}^s
n_i\Z v_i\subseteq\Lambda\,$$ of index coprime to $p$. Then for
$K={\bar \F}_p$ the isomorphism $j$ as in (\ref{j}) yields
bijections
$$j:\Char_{a-{\rm harm}} (\Lambda,K^\times)
\stackrel{\cong}{\longrightarrow} \Sigma_{a}\subseteq
(K^\times)^s\,$$ and
$$j:\Char_{a-{\rm harm}}
(\T_{\bar n},K^\times) \stackrel{\cong}{\longrightarrow}
\Sigma_a\cap\mu_{\bar n}\,,$$ where $\T_{\bar
n}:=\Lambda/\Lambda'$.
\ethm

\bcor\label{c1}
$$\dim\,{\rm Harm}_a \,(\T_{\bar n}, K)=\card\,\left(\Char_{a-{\rm
harm}}(\T_{\bar n}, K^\times)\right) \,$$$$=\card
\,(\Sigma_a\cap\mu_{\bar n})=\mult_{\lambda=0} \,(T_{a,
\Lambda'})\,.$$ \ecor

\bobs\label{o5} Reversing the logic we let $\Sigma$ be an
arbitrary affine algebraic subvariety  in the torus
$(K^\times)^s$. Thus $\Sigma$ can be defined by a finite sequence
$(p_j)$ of Laurent polynomials. When does $\Sigma$ possess a point
with a given torsion multi-order?

To answer this question, we pass to the associated system of
Laplacians $\Delta_{a_j},\,j=1,\ldots,t$, where $a_j:\Z^s\to K$
is the coefficient function
of the polynomial $p_j$. It is easily seen that the
orthogonal projection
$$\pi:\mathcal F (\Lambda, K)^{\Lambda'}\to \ker \,(\Delta_{\bar
a})\,$$ is given by
$$\pi=\prod_{j} \left(1-\Delta_{a_j}^q\right)\,$$
for a suitable $q=p^\alpha$.

If $t=1$ i.e., $\Sigma=\Sigma_a \subseteq (K^\times)^s$ is a
hypersurface, we can indicate a dynamical way to determine whether
$\Sigma_{a}\cap\mu_{\bar n}\neq\emptyset$.  Namely the latter
holds if and only if the following sequence of functions on the
quotient group $\T_{\bar n}=\Lambda/\Lambda'$  is not periodic:
$$f_0=\delta_0,\quad f_k=\Delta_a^k(\delta_0),\qquad k\in\N\,.$$
Indeed, in the latter case $f_{k+l}=f_k$ for certain
minimal $k,l$ with $k>0,l>0$, and
so $h=f_{k+l-1}-f_{k-1}$
is a nonzero harmonic function on $\T_{\bar n}$.
\eobs

\subsection{Winning boards and Artin's conjecture on primitive roots}
The following results were elaborated in Hunziker-Machiavelo-Park
\cite{HMP}. We formulate them in terms of existence of a nonzero
harmonic function on a toric square board with values in a Galois
field $\F_p$ as $p$ varies.

\bthm\label{5} {\rm (\cite{HMP})}
For an $n\times n$ torus $\T_{n,n}$
the following hold.

\smallskip

\noindent  (a) $$\forall n\ge 3\,\,\exists p\,:\, {\rm Harm}
\,(\T_{n,n},\F_p)\neq (0)\,.$$

\noindent (b) $${\rm Harm} \,(\T_{n,n},\F_p)\neq (0)\quad\forall p
\quad\Longleftrightarrow\quad n\equiv 0 \!\!\mod
3\quad\mbox{or}\quad n\equiv 0 \!\!\mod 5\,.$$

\smallskip

\noindent (c)  Except for at most 2 values of the prime $p$, the set
of all primes $l$ such that
$${\rm Harm} \,(\T_{l,l},\F_p)=(0)$$ is infinite.

\smallskip

\noindent (d) If $n=\frac{p\pm 1}{2}$ and $p\ge 23$ then ${\rm
Harm} \,(\T_{n,n},\F_p)\neq (0)$.

\smallskip

\noindent (e) We let $P_p$ denote the set
of all $n\in\N$ such that
${\rm Harm} \,(\T_{n,n},\F_p)\neq (0)$
 while for any proper divisor $d$ of $n$,
 ${\rm Harm} \,(\T_{d,d},\F_p)= (0)$.
Then $P_2$ and $P_3$ are infinite.
\ethm

\noindent The proof of (c) is based on a result of Heath-Brown
\cite{HB}, which concerns the following conjecture.

\medskip

\noindent {\bf Artin's conjecture on primitive roots} (1927; see
\cite{Mo,Mu}): {\it Every integer $n\neq -1$ which is not a square
is a primitive root modulo $l$ for an infinite number of primes
$l$.}

\medskip

\noindent Due to \cite{HB} this conjecture holds indeed for all
primes $n=p$ with at most 2 exceptions, and for all square-free
integers $n$ with at most 3 exceptions (see also \cite[\S 5]{Mo}).
For instance, at least one among the primes $2,3,5$ must satisfy
Artin's condition. However, no specific prime $p$ is known to
possess the Artin property.

\section{Harmonic functions on trees}
  We fix
 a field $K$ of characteristic $p>0$ and a finite graph $\Gamma$.
 According to Amin-Slater-Zhang \cite{ASZ}
 and Gravier-Mhalla-Tannier \cite{GMT},
 one simplifies $\Gamma$ by applying the following two surgeries:

\smallskip

\noindent $\bullet$ suppressing $p$ extremal vertices
$u_1,\ldots,u_p$ joint with a common neighbor $v$ together with
the edges $[u_i,v]$ ($i=1,\ldots,p$) as shown on the following
diagram:

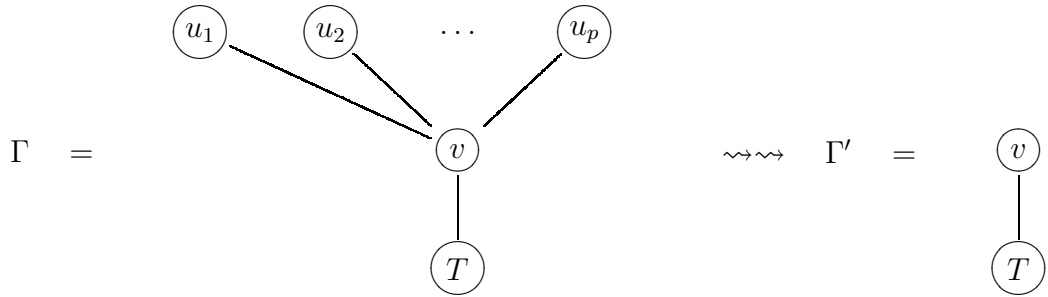
\begin{equation*}
\xymatrix{ & *++[o][F]{u_1} \ar@{-}@//[drr] & *++[o][F]{u_2}
\ar@{-}@//[dr] & \ldots
&  *++[o][F]{u_p} \ar@{-}@//[dl] & & &\\
\Gamma\quad= & & & *++[o][F]{v} \ar@{-}@//[d] &
 &\quad\rightsquigarrow\rightsquigarrow\quad  \Gamma'\quad= & *++[o][F]{v} \ar@{-}@//[d] \\
& & &*++[o][F]{T} & & & *++[o][F]{T} & &}
\end{equation*}

\vskip 0.3in

\noindent $\bullet$ suppressing an extremal linear branch of
length 2, say, $[u,v,w]$, where $u$ is an extremal vertex,
together with all edges joining $w$ with the rest $\G'$ of $\G$:

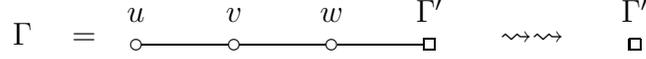
$$
\G\quad=\quad\co{u}\llin\co{v}\llin\co{w}\llin\boxo{\G'}\quad\quad
\rightsquigarrow\rightsquigarrow\quad \quad \boxo{\G'}\quad $$ The
following simple observation is essentially due to \cite{GMT}.

\bprop\label{p2} Let a graph $\Gamma'$ be obtained from $\Gamma$
by performing a surgery as above. Then any $K$-valued harmonic
function on $\G$ restricts to a $K$-valued harmonic function on
$\G'$ and vice versa, any $K$-valued harmonic function on $\G'$
extends uniquely to a $K$-valued harmonic function on $\G$. This
extension provides an isomorphism
$${\rm Harm}\,(\Gamma',K)\cong {\rm Harm}\,(\Gamma,K)\,.$$
\eprop

In particular, for a linear path $\PP_n$ with $n\ge 3$ vertices we
have
$${\rm Harm}\,(\PP_n,K)\cong {\rm Harm}\,(\PP_{n-3},K)\,.$$
Consequently, ${\rm Harm}\,(\PP_n,K)\neq (0)$ if and only if
$n\equiv 2 \mod 3$.

Every finite graph $\Gamma$ can be reduced, via a suitable
sequence of the surgeries as above, to a graph, say,
$\red_p\,(\Gamma)$ such that any extremal linear branch of
$\red_p\,(\Gamma)$ has length 1, and at each vertex of
$\red_p\,(\Gamma)$ there are at most $p-1$ such extremal linear
branches.

In case $p=2$ for a finite tree $\G$ such a graph
$\red_2\,(\Gamma)$ consists of isolated vertices, say,
$w_1,\ldots,w_t$ and isolated edges, say,
$[u_1,v_1],\ldots,[u_s,v_s]$. Any harmonic function $h$ on $\G$
takes value $0$ at $w_1,\ldots,w_t$, while $h(u_i)=h(v_i)$
$\forall i=1,\ldots,s$. Conversely, any such a binary function on
$\red_2\,(\Gamma)$ is harmonic and extends uniquely to a binary
harmonic function on $\G$. This leads to the following result.

\bcor\label{c3} For a finite tree $\G$, $\dim {\rm Harm}\,(\Gamma,
\F_2)$ is equal to the number $s$ of isolated edges in any
reduction $\red_2\,(\Gamma)$ of $\G$. \ecor

We wonder whether there exists an alternative (non-recursive)
combinatorial interpretation of the invariant $\dim {\rm
Harm}\,(\Gamma, \F_2)$ (it is related to the ``parity dimension"
as considered e.g., in \cite{ASZ}).

\brem\label{rem100} Let us note that the number of isolated
vertices in a reduction $\red_2\,(\Gamma)$ depends on the
reduction. Thus this is not in general an invariant of $\G$.
Indeed, a linear path with 3 vertices $\PP_3$ admits two different
reductions. One of them is empty, while the other one is $\PP_1$
and so consists of a single vertex. \erem

\bexa\label{e2} We say that a graph $\Gamma$ is odd if the degree
of each vertex of $\Gamma$ is. By virtue of Proposition 6 in
\cite{ASZ}, for any finite odd tree $\Gamma$
$$\dim {\rm Harm}\,(\Gamma, \F_2)=1\,.$$
The only nonzero binary harmonic function on $\Gamma$ is the
constant function $1$. However, $\dim {\rm Harm}\,(\Gamma,
\F_2)=\infty\,$ for any infinite locally finite odd tree
$\Gamma$.\eexa

\section{Appendix: Classical polynomials}
\subsection{Normalized Chebyshev polynomials and Fibonacci polynomials}
We recall that the Chebyshev polynomials of the first (second)
kind\footnote{With the conventions used e.g., in MAPLE.} satisfy
the following relations:
$$V_0=1,\quad V_1=x\quad\mbox{ and}\quad
V_n=2xV_{n-1}-V_{n-2}\quad\forall n\ge 2\,,$$ respectively,
$$U_0=1,\quad U_1=2x\quad\mbox{ and}\quad
U_n=2xU_{n-1}-U_{n-2}\quad\forall n\ge 2\,.$$ The normalized
Chebyshev polynomials of the first (second) kind $G_n\in\Z[x]$
($F_n\in\Z[x]$, respectively) are defined \cite{HMP} via
$$G_0(x)=2,\quad G_1(x)=x,\quad\mbox{and}\quad
G_n(x)=xG_{n-1}(x)-G_{n-2}(x) \quad\forall n\ge 2\,,$$
respectively, $$F_0(x)=0,\quad F_1(x)=1,\quad\mbox{and}\quad
F_n(x)=xF_{n-1}(x)-F_{n-2}(x)\qquad\forall n\ge 2\,.$$ The
Fibonacci polynomials $f_n\in\Z[x]$ are generated via the
recurrence relation:
$$f_0=0,\quad f_1=1,\quad f_{n}=xf_{n-1}+f_{n-2}\quad\forall n\ge 2\,.$$
They reduce to the Fibonacci numbers for $x=1$ and satisfy
identities similar to those known for the Fibonacci numbers. The
polynomials $F_n,G_n$ and $xf_n$ are even (odd) iff $n$ is. We
have $G_n(x)=2V_n\left(\frac{x}{2}\right)$ and
$$(-1)^n\det \left({\rm adj}\,(\T_n)-xI_n\right)
=C_n(x)=G_n(x)-2\qquad\forall n\ge 3\,,$$ respectively,
$$\det \left({\rm
adj}\,(\PP_{n-1}) -xI_{n-1}\right)=F_n(-x)\quad\forall n\ge 2\,.$$

\bprop\label{ch-hmp} \cite[\S 2]{Bi, HBJ, HMP} The normalized
Chebyshev polynomials $F_n$, $G_n$ and the Fibonacci polynomials
$f_n$ acquire the following properties:
\begin{enumerate}\item[$\bullet$]
$F_n(x+x^{-1})=\frac{x^n-x^{-n}}{x-x^{-1}}$ and
$G_n(x+x^{-1})=x^n+x^{-n}$. \item[$\bullet$]
$\gcd\,(F_m,F_n)=F_{\gcd\,(m,n)}$ and
$\gcd\,(f_m,f_n)=f_{\gcd\,(m,n)}$. \item[$\bullet$]
$F_m\,|\,F_n\quad\Longleftrightarrow\quad
f_m\,|\,f_n\quad\Longleftrightarrow\quad m\,|\,n$.
\item[$\bullet$] $F_{mn}=F_n\cdot (F_m\circ G_n)$ and
$G_{mn}=G_m\circ G_n$. \item[$\bullet$] $(x^2-4)F_m
F_n=G_{m+n}-G_{|m-n|}$ and $G_mG_n=G_{m+n}+G_{|m-n|}$.
\end{enumerate}\eprop

The next result deals with the irreducible factorization of the
Fibonacci polynomials.

\bprop\label{irr-fact} \cite[Corollary 2.3]{Le, JRS} There are
irreducible polynomials $\theta_n\in\Z[x]$ with nonnegative
coefficients, of degree $\deg \theta_n=\varphi(n)$ \footnote{Here
$\varphi$ stands as before for the Euler totient function.} such
that
$$f_n=\prod_{d|n}\theta_d\qquad\forall n\ge 1\,.$$ \eprop

\subsection{Dickson polynomials}
We recall \cite{LMT} that the Dickson polynomials
$D_n(x,a)\in\Z[x,a]$ and $E_n(x,a)\in\Z[x,a]$ of the first
(second) kind are defined recursively via:
$$D_0=2,\quad D_1=x,\quad\quad
D_{n+1}(x,a)=xD_n(x,a)-aD_{n-1}(x,a)\,,$$ and
$$E_0=1,\quad E_1=x,\quad\quad
E_{n+1}(x,a)=xE_n(x,a)-aE_{n-1}(x,a)\,,$$ respectively. They can
be also characterized by the identities:
$$ D_n(\mu_1+\mu_2,\mu_1\mu_2)=\mu^n_1+\mu^n_2\qquad \mbox{resp.,}
\qquad
E_n(\mu_1+\mu_2,\mu_1\mu_2)=\mu^{n+1}_1-\mu^{n+1}_2/(\mu_1-\mu_2)\,.$$
Furthermore, $E_{n-1}={D_n}'/n$. For $a=1$ the Dickson polynomials
specialize to the normalized Chebyshev polynomials:
$$G_n(x)=D_n(x,1)\qquad\mbox{and}\qquad F_{n+1}(x)=E_n(x,1)\,.$$
Similarly, $f_n(x)=\tilde f_n(x,1)$, where $\tilde
f(x,y)\in\Z[x,y]$ stands for the bivariate Fibonacci polynomials.
These are defined \cite{HL} by the recursion
$$\tilde f_0=0,\quad \tilde f_1=1,\quad \tilde f_{n}
=x\tilde f_{n-1}+y\tilde f_{n-2}\quad\forall n\ge 2\,.$$ The
polynomials $f_n$ and $\tilde f_n$ are irreducible over $\Q$ if
and only if $n$ is prime \cite{HL, WP}. An analog of Proposition
\ref{irr-fact} also holds for $\tilde f_n$ \cite{JRS}.

\subsection{Reduction to a positive characteristic}
Given a prime $p$ and a polynomial $F\in \Z[x]$, we let
$F^{(p)}\in\F_p[x]$ denote the reduction  of $F$ modulo $p$. The
Dickson polynomials reduced modulo $p$ satisfy the relations
\cite{BZ}:
$$D^{(p)}_{p^\alpha m}= {(D^{(p)}_m)}^{p^\alpha}
\qquad\mbox{resp.}\qquad E^{(p)}_{p^\alpha m-1}=
{(E^{(p)}_{m-1})}^{p^\alpha}(x^2-4a)^{\frac{p^\alpha-1}{2}}\,,$$
where $m\not\equiv 0\mod p$. Similarly, for the reduction
$F_m^{(p)}$ of the normalized Chebyshev polynomials of the second
kind we have

\bprop\label{ch-hmpp} \cite[\S
2]{HMP}\begin{enumerate}\item[$\bullet$]
$F_m^{(p)}\,|\,F_n^{(p)}\quad\Longleftrightarrow\quad m\,|\,n$.
\item[$\bullet$] $F_{p^km}^{(p)}=F_{p^k}^{(p)}\cdot
{(F_m^{(p)})}^{p^k}$, where $F_{p^k}^{(p)}=(x^2-4)^{(p^k-1)/2}$.
\item[$\bullet$]
$F_{(p^k-1)/2}^{(p)}F_{(p^k+1)/2}^{(p)}=\frac{x^{p^k}-x}{x^2-4}$
if $p\neq 2$.
\end{enumerate}\eprop

\bcor\label{ch-irr} \cite[\S 2]{HMP} Every irreducible polynomial
$\tau\in \F_p[x]$ of degree $k$ occurs as a factor of
$F_{p^k-1}^{(p)}F_{p^k+1}^{(p)}$ if $p=2$ or of
$F_{(p^k-1)/2}^{(p)}F_{(p^k+1)/2}^{(p)}$ if $p\neq 2$. \ecor

Reducing the polynomials  $G_n$ modulo $2$ yields the
Chebyshev-Dickson polynomials $T_n\in\F_2[x]$. Actually
$T_n=C_n^{(2)}=G_n^{(2)}=xF_n^{(2)}=xf_n^{(2)}$ for $n\ge 3$
\footnote{Although $C_0\equiv T_0,\,\,C_1\equiv T_1  \mod 2$,
however $C_2=x^2-1\not\equiv T_2=x^2 \mod 2$.}. They can also be
defined recursively:
$$T_0=0,\qquad T_1=x,\qquad T_{n+1}=xT_n+T_{n-1}\,,$$
or, alternatively, via the relation $$T_n(x)\equiv
xU_{n-1}\left(\frac{x}{2}\right) \mod 2\,,$$ where $U_n\in\Z[x]$
stands for the $n$th Chebyshev polynomial of the second kind.

\bprop\label{ch-hmp} \cite{GKW, SB, Su} The Chebyshev-Dickson
polynomials $T_n\in \F_2[x]$ acquire the following properties:
\begin{enumerate}
\item[$\bullet$]  $T_n(x+x^{-1})=x^n+x^{-n}$. \item[$\bullet$]
$T_m\circ T_n=T_{mn}$. \item[$\bullet$]
$\gcd\,(T_m,T_n)=T_{\gcd\,(m,n)}$. \item[$\bullet$]
$T_m\,|\,T_n\quad\Longleftrightarrow\quad m\,|\,n$.
\item[$\bullet$]  $T_{2^km}=T_{m}^{2^k}$. \item[$\bullet$]
$T_{2^k-1}T_{2^k+1}=(x^{2^k-1}-1)^2$.
\end{enumerate}\eprop

\end{document}